\documentclass[11pt,a4paper]{article}

\usepackage{amsthm,amssymb,hyperref}
\usepackage[all]{xy}
\title{A Fourier-Mukai approach \\ to spectral data for instantons}
\author{Marcos Jardim \\ University of Massachusetts at Amherst \\
Department of Mathematics and Statistics \\ Lederle Graduate Research Tower \\
Amherst, MA 01003-9305 USA \\ \\ Anthony Maciocia \\
Department of Mathematics and Statistics\\University of Edinburgh \\
James Clark Maxwell Building\\King's Buildings\\ Mayfield Road\\
Edinburgh, EH9 3JZ, UK }

\newcommand{\seta}{\rightarrow}\newcommand{\lseta}{\longrightarrow}
\newcommand{\torus}{\mathbb{T}} \newcommand{\dual}{\hat{\mathbb{T}}}
\newcommand{\dv}{\hat{V}} \newcommand{\dw}{\hat{W}}
\newcommand{\ext}{\mathop{\rm Ext}\nolimits} \newcommand{\pp}{\mathcal{P}}
\newcommand{\Hom}{\mathop{\rm Hom}\nolimits}
 \newcommand{\real}{\mathbb{R}}
\newcommand{\cpx}{\mathbb{C}} 
\newcommand{\proj}{\mathbb{P}^1} 
\newcommand{\cF}{\mathcal{F}} 
\newcommand{\cM}{\mathcal{M}} 
\newcommand{\oo}{\mathcal{O}} \newcommand{\cS}{\mathcal{S}}
\newcommand{\nn}{{\cal N}}
\newcommand{\cL}{\mathcal{L}} \newcommand{\hE}{\hat{E}}

\newcommand{\projb}{\mathcal{B}}
\newcommand{\wit}[1]{-WIT$_{#1}${}}
\newcommand{\po}{P} 
\newcommand{\ch}{\mathop{\rm ch}\nolimits} 

\newtheorem{thm}{Theorem} \newtheorem{lem}[thm]{Lemma}
\newtheorem{prop}[thm]{Proposition} \newtheorem{cor}[thm]{Corollary}
\theoremstyle{remark} \newtheorem*{rem}{{\bf Remark}} 
\newtheorem*{dfn}{{\bf Definition}} 

\begin{document}

\maketitle

\begin{abstract}
We study $SU(r)$ instantons on elliptic surfaces with
a section and show that they are in one-one correspondence with
spectral data consisting of a curve in the dual elliptic surface and a
line bundle on that curve. We use relative Fourier-Mukai transforms to
analyse their properties and, in the case of the K3 and abelian
surface, we show that the moduli space of instantons has a natural
Lagrangian fibration structure with respect to the canonical complex
symplectic structures.
\end{abstract}

\newpage

\section{Introduction} \label{intro}

The mathematical study of gauge theory was born some 25 years ago, and one of
the first main results is just as old, namely the ADHM construction of
instantons over the Euclidean 4-space. Since then, several different types
of Euclidean instantons have been studied from various points of view:
monopoles, calorons, Higgs pairs, doubly-periodic instantons, to name a few.
The common feature in the study of such objects is the so-called {\em Nahm
transform}.  It relates instantons on $\real^4$ which are invariant under a
subgroup of translations  $\Lambda\subset\real^4$ to instantons on the dual
$(\real^4)^*$ which are invariant under the dual subgroup of translations
$\Lambda^*$ (see \cite[Section 7]{Nk} for a detailed exposition). The
Nahm transform has many interesting properties, but perhaps its greatest
virtue is that it often converts the difficult problem of solving nonlinear
PDE's into a simpler problem involving only ODE's or just vector spaces
and linear maps between them.

More recently, the so-called {\em Fourier-Mukai transforms} have generalised the
Nahm transform, bringing it to the realms of algebraic geometry and
derived categories. In this paper, we shall study a particular version
of such transforms. It is defined for torsion-free sheaves $E$
over relatively minimal elliptic surfaces $X\stackrel{\pi}{\seta}B$
with a section, transforming them into torsion sheaves on
$J_X\stackrel{\hat{\pi}}{\seta}B$, the relative Jacobian of $X$, which
are supported over the spectral curves studied by
Friedman-Morgan-Witten \cite{FMW1,FMW2} for physical reasons. These
constructions are briefly reviewed in Sections \ref{spec} and
\ref{fm}, and we observe that there is a correspondence between
instantons on a 4-torus, spectral curves with line bundles on it,
and instantons on the dual 4-torus.

The bulk of the paper is contained in Section \ref{preserve}, where we compare
the $\mu$-stability of $E$ with the concept of stability introduced by Simpson
\cite{S} for its transform. In particular, we show that {\em $E$ is $\mu$-stable
and locally-free if and only if its transform is stable in the sense of Simpson.}

As an application of these ideas, we shed new light into the moduli space of
irreducible $SU(r)$-instantons over elliptic K3 and abelian surfaces, showing
in Sections \ref{fibre} and \ref{hyper} that it has the structure of a complex
Lagrangian fibration. The case of rank 2 instantons on an elliptic surface is
treated in Friedman's book \cite{F}. We extend the results to
higher rank and provide some more details about the structure of the
fibration. The proofs are greatly simplified using the Fourier-Mukai
technology and this enables us to go further in our description of the
moduli spaces.

\paragraph{Notation.}
The elliptic surface $X$ will be polarised by $\ell$ which can be
any polarisation (for the product of elliptic curves it is convenient
to choose the sum of the elliptic curves).
If $F_x$ denotes the fibre of $X\stackrel{\pi}{\seta}B$ passing
through $x\in X$, then  $\hat{\pi}^{-1}(\pi(x))=\hat{F}_{\pi(x)}$ is
the dual of $F_x$. Let $J_X$ denote the dual elliptic surface
parametrising flat line bundles on the fibres of $\pi$. It also
fibres over $B$ and has a canonical section given by the trivial
bundles. We shall denote its fibres by $\hat F_b$, for $b\in B$.

Given a sheaf $E$ on $X$, we write its Chern character as a triple:
$$ \ch(E)=({\rm rank}(E),c_1(E),\ch_2(E)), $$
where $\ch_2(E)=\frac{1}{2}c_1(E)^2-c_2(E)$.

By $\pp\seta X\times_B J_X$ we mean the
relative Poincar\'e sheaf parametrising flat line bundles on the
smooth fibres. Given $p\in J_X$, let $\iota_p$ be the
inclusion of $F_{\hat{\pi}(p)}$ inside $X$. We define:
$$ E(p) = \iota_p^*(E\otimes\pp_p) $$
where $\pp_p$ denotes the torsion sheaf on $X$ corresponding to $p\in J_X$.

Finally, $\cM_X(r,k)$ (or just $\cM$) will denote the moduli space of
$\mu$-stable locally-free sheaves over $X$, with Chern classes $(r,0,-k)$,
where $r,k>1$. Note that this coincides with the moduli space of $SU(r)$
instantons of charge $k$.


\section{Spectral Data for instantons} \label{spec}

Let $\cal E$ be a vector bundle over the elliptic surface $X$, and
let $A$ be an anti-self-dual $SU(r)$ connection on $\cal E$. Then $A$ induces
a holomorphic structure $\overline{\partial}_A$ on $\cal E$; we denote
by $E$ the associated holomorphic vector bundle.

Given a holomorphic vector bundle $E\seta X$ as above, we say that $E$
is {\em generically fibrewise semistable} if its restriction to
a generic elliptic fibre is semistable; $E$ is said to be
{\em fibrewise semistable} if its restriction to every elliptic
fibre is semistable. We shall also say that $E$ is {\em regular}
if $h^0(F_{\hat{\pi}(p)},E(p))\leq 1$ for all $p\in J_X$. Observe
that regularity is a generic condition.

Now assume that $E$ is a regular holomorphic vector bundle over $X$;
we define the {\em instanton spectral curve} with respect to the projection
$X\stackrel{\pi}{\seta}B$ in the following manner:
\begin{equation} \label{spec.dfn}
S = \{ p\in J_X\ |\ h^1(F_{\hat{\pi}(p)},E(p))\neq0 \}
\end{equation}
It is not difficult to see that $S\stackrel{\pi}{\seta}B$ is a
branched $r$-fold covering map.

To define the second part of the instanton spectral data, recall that
$\chi(E(p))=0$ for every  $p\in J_X$. We define a line bundle
$\cL$ on $S$ by attaching the vector space $H^1(F_{\hat{\pi}(p)},E(p))$
to the point $p\in S$. Alternatively, consider the diagram:
\begin{equation} \label{maps1} \xymatrix{
X\times_B J_X& X\times_B S\ar[l]_-{\tau_1} \ar[r]^-\sigma \ar[d]^-{\tau_2} & S \\ & X
} \end{equation}
where $\tau_1$ is the inclusion map and $\tau_2,\sigma$ are the
obvious projections.
We define $\cL$ as follows:
\begin{equation} \label{lbdefn1}
\cL \ = \ R^1\sigma_*(\tau_2^*E\otimes\tau_1^*\pp)
\end{equation}

Geometrically, it is not difficult to see that $\cL$ can be regarded as a
{\em bundle of cokernels} associated to a family of coupled Dirac operators
parametrised by $S$. For the simple case $X=T\times\proj$, see \cite{J3}.
This follows from the natural identification between $H^1(F_{\hat{\pi}(p)},E(p))$
and the cokernel of the Dolbeault operator $\overline{\partial}_A|_{F_{\hat{\pi}(p)}}$.

As we will see in Section \ref{fm} below, the above construction
is invertible, and the original holomorphic bundle $E$ can be
reconstructed from its spectral pair $(S,\cL)$. However, the
spectral pair contains only the holomorphic information, so some
extra information is needed in order to reconstruct the original
instanton connection. One possibility is to construct a connection
on $\cL\seta S$, as it was done in \cite{J3}; see also \cite{LYZ}
for a related problem. In this paper we will follow a different
path, which will be described in Section \ref{preserve} below.


\section{Spectral Data and the Relative \\ Fourier-Mukai Transform} \label{fm}

Given a sheaf $E$ on the relatively minimal elliptic surface $X\rightarrow B$
with a section, we define its relative Fourier-Mukai transform to be the
{\em complex} of sheaves $\Phi(E)$ on $J_X$ given by:
\begin{equation} \label{phi}
\Phi(E) = R\hat{\pi}_* \big( \pi^*E \otimes \pp \big)
\end{equation}
where $X\stackrel\pi\longleftarrow X\times_B
J_X\stackrel{\hat\pi}\longrightarrow J_X$ are the projection maps and
$\pp$ is the relative Poincar\'e sheaf as before.

We say that $E$ is $\Phi$-WIT$_n$ if $\Phi^i(E)=0$ unless $i=n$, where
$\Phi^i(E)$ denotes the $i$th homology of the complex representing
$\Phi(E)$ which is well defined up to isomorphism.

It can be shown that (\ref{phi}) gives an equivalence between $D(X)$ and
$D(J_X)$, the derived categories of bounded complexes of coherent sheaves.
Its inverse functor is given by:
\begin{equation} \label{invphi}
\hat{\Phi}(L) = R\pi_* \big( \hat{\pi}^*L \otimes \pp^\vee \big)
\end{equation}
so that $\hat{\Phi}^0(\Phi^1(E))=E$ and $\Phi^1(\hat{\Phi}^0(L))=L$

This is one of a series of Fourier-Mukai transforms for elliptic
surfaces parameterized by SL$_2(\mathbb{Z})$. For a description of
these and their invertibility, see \cite{B}. The rank and fibre
degree of the transform $\Phi(E)$ is related to the rank and fibre
degree of $E$ by multiplication by the element of
SL$_2(\mathbb{Z})$. The case we shall be interested in this paper
is given by the matrix
\[\left(\begin{array}{cc}0&1\\-1&0\end{array}\right)\in\mathrm{SL}_2(\mathbb{Z})\]

Let us now recall the details of the definition of stability
(due to Mumford and Takemoto) of torsion-free sheaves:

\begin{dfn}
A torsion-free sheaf $E$ over a polarised surface $(X,\ell)$ is said
to be $\mu$-\textit{stable} (respectively, $\mu$-\textit{semistable})
if for all proper subsheaves $F$ of $E$, $\mu(F)<\mu(E)$ (respectively,
$\mu(F)\leq\mu(E)$), where $\mu(\bullet)=c_1(\bullet)\cdot\ell/r(\bullet)$
is called the \textit{slope}.
\end{dfn}

Recall that since we are working over surfaces, we can assume that
the destabilizing subsheaf is locally-free and that the quotient
$E/F$ is torsion-free. Moreover, we also may choose $F$ to be
$\mu$-stable.

If $E$ is $\mu$-stable and locally-free, that is if $E$ is a holomorphic
vector bundle arising from an irreducible instanton connection via the
Hitchin-Kobayashi correspondence, it is easy to see that $E$ is $\Phi$\wit1.
Hence $\Phi(E)=\Phi^1(E)$ is not a complex of sheaves, but simply a sheaf on
$J_X$. This also applies if $E$ is not necessarily locally-free but is still
torsion-free, as we shall see this in Section \ref{preserve}.

The following Proposition brings together the instanton spectral data described
in the previous Section with the Fourier-Mukai transform.

\begin{prop}
Let $\cal E$ be a vector bundle with an irreducible $SU(r)$ instanton
connection over an elliptic surface $X$. Let $E$ be the associated
$\mu$-stable, regular holomorphic vector bundle with $c_1(E)=0$.
Then $\Phi^1(E)$ is a torsion sheaf supported over the instanton
spectral curve $S\subset J_X$. Moreover, the restriction of $\Phi^1(E)$
to its support is naturally isomorphic to $\cL$.
\end{prop}
\begin{proof}
Let $j:S\to J_X$ be the inclusion map. This fits into a commuting
diagram of maps introduced in (\ref{maps1}):
$$ \xymatrix{
X\times_B S\ar[r]^-{\tau_1}\ar[d]^-{\sigma} &
X\times_B J_X\ar[d]^-{\hat\pi} \\
S\ar[r]^-j & J_X } $$
Then
\begin{eqnarray*}
j_*\cL&\cong&R(j_*\sigma_*)(\tau_2^*E\otimes\tau_1^*\pp)\\
&\cong&R\hat\pi_*R\tau_{1*}(\tau_2^*E\otimes\tau_1^*\pp)\\
&\cong&R\hat\pi_*(\pp\otimes R\tau_{1*}\tau_2^*E).
\end{eqnarray*}
But $\tau_2=\pi\tau_1$ (see diagram \ref{maps1}), so
$R\tau_{1*}\tau_2^*E\cong\pi^*E\otimes\oo_{X\times_B S}$. Hence,
we have a natural map $\Phi^1(E)\to j_*\cL$. Since the fibres of
$\Phi^1(E)$ and $j_*\cL$ are naturally isomorphic, this must be an
isomorphism. Hence, $j^*\Phi^1(E)=j^*j_*\cL=\cL$ as required.
\end{proof}


\section{Instantons on flat 4-tori} \label{product}

Let $\torus$ be a flat 4-dimensional torus, with a fixed complex
structure. Denote by $\dual$ the corresponding dual torus, which
inherits a flat metric and a complex structure from $\torus$ polarized
by the associated K\"ahler class.

We start by considering an irreducible $SU(r)$ bundle ${\cal E}\seta\torus$ with 
$c_1({\cal E})=0$ and $c_2({\cal E})=k$, plus an instanton connection $A$.
Equivalently, we can adopt an algebraic geometric point of view and look at
the associated $\mu$-stable holomorphic vector bundle $E\seta\torus$ of rank
$r$ and the same Chern classes. Stability now implies the irreducibility
hypothesis, which in turns implies that $h^0(\torus,E)=h^2(\torus,E)=0$,
so that $h^1(\torus,E)=k$.

The {\em Nahm transformed bundle} $\hE\seta\dual$ is then constructed as
follows. Let $\po\seta\torus\times\dual$ be the Poincar\'e line bundle,
and consider the natural projection maps:
\begin{equation} \label{nt}
\torus\ \stackrel{p}{\longleftarrow}\ \torus\times\dual\
\stackrel{\hat{p}}{\longrightarrow}\ \dual
\end{equation}
Then $\hE=R^1\hat{p}_*(p^*E\otimes\po^\vee)$. The fibre of $\hE$ over
$\xi\in\dual$ can be canonically identified with $\hE_\xi=H^1(\torus,E\otimes L_\xi)$,
where $L_\xi\seta\torus$ denotes the flat holomorphic line bundle associated with
$\xi\in\dual$.

Conversely, it can be shown that $E=R^1p_*(\hat{p}^*\hE\otimes\po)$, with fibres
canonically identified with $E_z=H^1(\dual,\hE\otimes L_z)$, where
$L_z\seta\dual$ denotes the flat holomorphic bundle associated with
$z\in\torus$.

\begin{thm}
$\hE$ is a $\mu$-stable holomorphic vector bundle of
rank $k$, flat determinant and $c_2(\hE)=r$. Therefore, since
it is invertible, the Nahm transform is a bijective correspondence
$\cM_\torus(r,k)\seta\cM_{\dual}(k,r)$.
\end{thm}

For the proof, we refer to \cite{DK}. The algebro-geometric version is given in
\cite{Mac}. In particular, in the $k=1$ case, the Theorem tells us that there
are no stable holomorphic bundles with this Chern character, thus precluding
the existence of instantons with unit charge. For larger $k$, this result is not
really helpful as it only converts the problem of constructing instantons/stable
bundles over $\torus$ into the same problem over $\dual$. A better understanding
of the moduli space of $SU(r)$ instantons on $\torus$ with $k>1$ requires a new tool,
namely the Relative Fourier-Mukai transform introduced in the previous Section.

Now assume that $\torus=V\times W$, i.e. $\torus$ is given by the product of two
elliptic curves $V$ and $W$. We shall denote by $\dv$ and $\dw$ the respective
Jacobian curves, so that  the dual torus is given by $\dual=\dv\times\dw$.
This case is particularly interesting because we have
two fibration structures leading apparently to two distinct spectral curves, and we
would like to understand the relation between them.

Regarding $\torus$ as an elliptic surface over $V$, consider the diagram:
\begin{equation} \label{fibv} \xymatrix{
V\times W & V\times W\times \dw \ar[l]_-{\pi_V} \ar[r]^-{\hat{\pi}_V}
& V\times \dw
} \end{equation}
where $\pi_V$ is the projection onto the first and second factors, while
$\hat{\pi}_V$ is the projection onto the first and third factors. Let also
${\bf P}_W\seta W\times\dw$ be the Poincar\'e line bundle. Then
\begin{equation} \label{lbdefn2}
\cL \ = \ R^1\hat{\pi}_{V*}(\pi_V^*E \otimes {\bf P}_W)
\end{equation}
is a torsion sheaf on $V\times\dw$, supported over the instanton
spectral curve $S$.

On the other hand, regard $\dual$ as an elliptic surface over $\dw$ and
consider the Nahm transformed bundle $\hE\seta\dual$. Based on the diagram:
\begin{equation} \label{fibw} \xymatrix{
\dv\times\dw & V\times V\times \dv\times \dw \ar[l]_-{\pi_W} \ar[r]^-{\hat{\pi}_W}
& V\times \dw
} \end{equation}
where $\pi_W$ is the projection onto the second and third factors, while
$\hat{\pi}_W$ is the projection onto the first and third factors. Let also
${\bf P}_V\seta V\times\dv$ be the Poincar\'e line bundle. Then
\begin{equation} \label{lbdefn3}
\nn  \ = \ R^1\hat{\pi}_{W*}(\pi_W^*E \otimes {\bf P}_V)
\end{equation}
is a torsion sheaf on $V\times\dw$, supported over the {\em dual} instanton
spectral curve $R$.

Since the bundles $E$ and $\hE$ are related via Nahm transform, it seems natural
to ask how are the two sets of spectral data $(S,\cL)$ and $(R,\nn)$ are related
to one another.

Regarding $\torus$ as an elliptic surface over $V$, we have just seen that the
spectral data $(S,\cL)$ is encoded into the torsion sheaf $\Phi^1(E)$ on
$V\times\dw$. On the other hand, regarding $\dual$ as an elliptic surface over
$\dw$ let us denote the Fourier-Mukai transform of $\hE$ by $\Psi(\hE)$. It is
a torsion sheaf on $V\times\dw$ encoding the dual spectral data $(R,\nn)$. Let also
$\hat{\Psi}$ denote the inverse of $\Psi$.

We are finally in a position to prove our first main result:

\begin{thm} \label{match}
The spectral pairs $(S,\cL)$ and $(R,\nn)$ are equivalent, in the sense that
$\cL$ and $\nn$ can be canonically identified as sheaves on $J_X$.
\end{thm}

\begin{proof}
If we define the functor $\cF:D(\torus)\seta D(\dual)$ as follows:
$$ \cF(E) = R \hat{p}_*(p^*E\otimes \po^\vee) $$
where $p$ and $\hat{p}$ are the projection maps in (\ref{nt}), then
$\hE=\cF^1(E)$. The functor $\cF$ is Mukai's original Fourier
transform introduced in \cite{M1}. Using this notation, it is enough to show
that $\hat{\Psi}\circ\Phi=\cF$.

To see this it suffices to show that
$\hat\Psi(\Phi(\oo_x))=\hat{\Psi}(\pp_x)=\po_x$ over $\dual$. But
$\pp_x\otimes\pp^\vee_y$ is supported at the intersection of
translates of $V$ and $\hat W$. Hence, $\pp_x$ is $\Psi$-WIT$_0$
and $\hat\Psi(\pp_x)$ is a flat line bundle. On the other hand,
the properties of $\Psi$ and $\Phi$ imply that translating $x$ to
$y$ twists $\hat{\Psi}(\pp_x)$ by $\po_{y-x}$ and so by
normalizing the Poincar\'e bundles appropriately, we have
$\hat{\Psi}(\pp_x)=\po_x$.
\end{proof}

Similarly, one could regard $\torus$ as an elliptic surface over $W$
and $\dual$ as an elliptic surface over $\dv$. The analogue of Theorem
\ref{match} would again hold, so that the corresponding spectral data
would also coincide as torsion sheaves in $\dv\times W$.

More generally, one has the following commuting diagram of derived
categories:

\begin{equation} \label{commdiag} \xymatrix{
D(V\times W) \ar[r]^\Phi \ar[d]^\Upsilon \ar[dr]^\cF &
D(V\times\dw) \ar[d]^\Psi \\
D(\dv\times W) \ar[r]^\Xi & D(\dv\times\dw)
} \end{equation}

Of course, the functor $\Upsilon$ is given by the Fourier-Mukai transform
regarding $\torus$ as an elliptic surface over $W$, while the functor $\Xi$
is the Fourier-Mukai transform regarding $\dual$ as an elliptic surface
over $\dv$.

\begin{rem}
One concludes from Theorem \ref{match} that there is an
equivalence between the following three objects: instantons on a
flat 4-torus (i.e. $\mu$-stable bundles over $\torus$), spectral
data, and instantons on the dual 4-torus (i.e. $\mu$-stable
bundles over $\dual$). Such ``circle of ideas" has been previously
establish for monopoles \cite{H1}, doubly-periodic instantons
\cite{J3} and periodic monopoles \cite{CK}. Indeed, one expects
that a similar scheme will hold for all translation  invariant
instantons on $\real^4$. \vskip18pt 
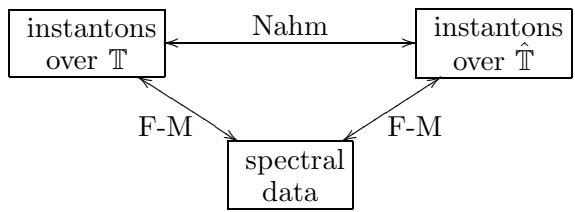
\centerline{ \xymatrix{
*+[F]\txt{ instantons \\ over $\torus$ }
\ar@{<->}[rr]^{\txt{Nahm}} \ar@{<->}[dr]_{\txt{F-M}}
& & *+[F]\txt{ instantons \\ over $\dual$ } \ar@{<->}[dl]^{\txt{F-M}} \\
& *+[F]\txt{ spectral \\ data } & }}
\vskip18pt
\end{rem}


\section{Preservation of Stability} \label{preserve}

We now aim at establishing the link between the stability
of $E$ and the stability of $L=\Phi^1(E)$. For torsion sheaves,
we have a concept of stability due to Simpson \cite{S}:

\begin{dfn}
A torsion sheaf $L$ on a polarised variety $(X,\ell)$ is
said to be $p$-\textit{stable} (respectively, $p$-\textit{semistable})
if it is of pure dimension (i.e. the support of all subsheaves have the
same dimension) and for all proper subsheaves $M$ of $L$ we have
$p(M,n)<p(L,n)$ (respectively, $p(M,n)\leq p(L,n)$) for all
sufficiently large $n$, where $p(L,n)$ denotes the reduced Hilbert
polynomial of $L\otimes\oo(n\ell)$.
\end{dfn}

In particular, a torsion sheaf which is supported on an irreducible
curve and whose restriction has rank one is automatically $p$-stable
with respect to any polarisation.

Now let $X\seta B$ be a (relatively minimal) elliptic surface
with a section $\sigma$; let $E$ be a torsion-free  sheaf on $X$ with
$\ch(E)=(r,0,-k)$.

\begin{prop} \label{one}
Suppose that the restriction of $E$ to some smooth fibre $F_s$ is
semistable. Then $E$ is $\Phi$\wit1 and $L=\Phi^1(E)$ has Chern character
$(0,k\hat{f}+r\hat{\sigma},0)$, where $\hat{\sigma}$ is the
class given by the zero section of $J_X$. In particular, $p(L,n)=n$.
\end{prop}
\begin{proof}
First note that if $E$ is torsion-free then
$\Phi^0(E)=\hat{\pi}_* \big( \pi^*E \otimes \pp \big)$
is also torsion-free. Hence its support consists is either empty
or consists of the whole relative jacobian surface $J_X$.
Moreover, if $E$ is not $\Phi$\wit1 then the support of
$\Phi^0(E)$ is contained in the support of $\Phi^1(E)$ since
$\chi(E|_{F_x})=0$ for all $x\in X$.

However, $\Phi^1(E)$ cannot be
supported on the whole surface since it does not contain the
whole fibre $\hat{F}_s$ where $E|_{F_s}$ is semistable. Thus
$\Phi^0(E)=0$ and $E$ is $\Phi$\wit1.

Now $\Phi^1(\oo_X)$ is the trivial line bundle supported on the zero section
of $J_X$ and $\Phi^0(\oo_x)$ is a flat line bundle supported on (a
divisor in) $\hat f$. Hence, since $\ch(E)=r\cdot\ch(\oo_X) + k\cdot\ch(\oo_x)$,
we conclude that:
$$ \ch(\Phi^1(E)) = r\cdot\ch(\Phi^1(\oo_X)) +
   k\cdot\ch(\Phi^0(\oo_x)) = (0,k\hat{f}+r\hat{\sigma},0) $$
as required.
\end{proof}

Note that the rank and fibre degrees of $L$ are given by $r(L)=c_1(E)\cdot f=0$
and $c_1(L)\cdot\hat f=r$ from the definition of the relative transform.

Let $S$ denote the support of $\Phi^1(E)$; the statement above does not imply
that $S$ is a divisor in $J_X$, since it might contain some 0-dimensional
components (i.e. isolated points).

\begin{lem} \label{two}
If $E$ is generically fibrewise semistable, then $S={\rm supp}\Phi^1(E)$
has no 0-dimensional components. Moreover, $L=\Phi^1(E)$ is of pure
dimension one.
\end{lem}
\begin{proof}
Suppose that $q\in S$ is an isolated point. It cannot belong to a fibre $\hat{F}_x$
such that $E|_{F_x}$ is semistable, since the restriction to every nearby fibre
is also semistable and the restriction of $E$ to the fibres varies
holomorphically.

So $E|_{F_q}$ is unstable. But $S$ contains all such fibres, since
$h^1(F_{\hat{\pi}(p)},E(p))>0$ for all $p\in\hat{\pi}^{-1}(\pi(q))$.

To obtain the second statement, we must show that $L$ has no proper
subsheaves supported on points. For a contradiction, suppose it does,
and let $F$ be a subsheaf of $L$ with 0-dimensional support. Then
$\hat{\Phi}^0(F)$ is a torsion subsheaf of $\hat{\Phi}^0(L)=E$,
thus contradicting the hypothesis that $E$ is torsion-free.
\end{proof}

Thus, we will call $S$ the {\em spectral curve} associated to the
generically fibrewise semistable torsion-free sheaf $E$, since it
generalizes our previous definition for regular locally-free
sheaves. As before, let us denote by $\cL$ the restriction of
$\Phi^1(E)$ to its support; in general, it is a coherent sheaf on
$S$. Note that ${\rm deg}(\cL)=0$.

\begin{rem}
If the restriction of $E$ to some smooth fibre has no sections, then
semicontinuity implies that $E$ is generically fibrewise semistable,
and hence $\Phi$\wit1. It is somewhat surprising that the spectral
data can actually be defined in an interesting, meaningful way under
such mild conditions.
\end{rem}

Summing up, we also conclude:

\begin{cor}
If $E$ is a torsion-free and $\Phi$\wit1, then $\Phi(E)=\Phi^1(E)$ is a
torsion sheaf of pure dimension one.
\end{cor}

Conversely, we also have:
\begin{lem} \label{L}
If $L$ is a $\hat{\Phi}$\wit0 torsion sheaf of pure dimension one
on $J_X$ with $ch(L)=(0,k\hat{f}+r\hat{\sigma},0)$, then
$\hat{\Phi}^0(L)$ is torsion-free.
\end{lem}

\begin{proof}
Suppose that the support of $L$ decomposes as $\Sigma+F$, where
$F$ is the sum of fibres.
When this happens, we have a subsheaf $K$ of $L$ supported on
$F$ with degree 0. Assume that $F$ is the maximal such effective
subdivisor of $S$. Let $Q=L/K$. Then $Q$ is $\hat\Phi$\wit0 and $K$ is
$\hat\Phi$\wit1. The resulting sequence after transforming with
$\hat\Phi$ is then
\begin{equation} \label{sqcL}
0\to E\to E^{**}\to\oo_Z\to 0
\end{equation}
where $E=\hat{\Phi}^0(L)$, for some zero dimensional subscheme $Z$.
\end{proof}

\begin{rem}
Notice that applying $\Phi$ to the sequence (\ref{sqcL})
shows that if $E$ is torsion-free but not locally-free,
then the support of $L$ contains fibres.
\end{rem}

The following lemma due to Bridgeland \cite{B} characterizes
$\hat{\Phi}$\wit0 sheaves on $J_X$:

\begin{lem} \label{tom's}
A sheaf $L$ on $J_X$ is $\hat{\Phi}$\wit0 if and only if
$$ {\rm Hom}(L,\pp_x)=0 \ \ \ \forall x\in X $$
\end{lem}

\subsection{Suitable polarizations}

Now let $\ell$ be a polarisation of the elliptic surface $X$, and let
$\hat{\ell}$ be the induced polarisation on $J_X$. If $\ell$ is
arbitrary, it is not difficult to see that a $\mu$-unstable torsion-free
sheaf on $X$ can have a $p$-stable transform.

Indeed, let $X$ be an elliptic surface whose fibres are all smooth, and
let $E$ be the locally-free sheaf given by the following extension:
$$ 0 \seta \oo(-\sigma+d f) \seta E \seta \oo(\sigma-d f) \seta 0 $$
Clearly, $c_1(E)=0$ and $c_2(E)=-\sigma^2+2d$. For $d$ sufficiently
large, $c_2(E)>0$ and $\ell\bullet(-\sigma+d f)=
-(\ell\bullet\sigma)+ d(\ell\bullet f)>0$, so $E$ is
$\mu$-unstable with respect to $\ell$.

On the other hand, notice that the restriction of $E$ to each fibre
is an extension of $Q\in{\rm Pic}^d(T)$ by its dual. Thus, the bundles
obtained by the above extension are generically regular. This means
$\Phi^1(E)$ is supported on a smooth, irreducible curve, so that
$\Phi^1(E)$ is necessarily $p$-stable with respect to $\hat{\ell}$.

Therefore, we can only expect the Fourier-Mukai transform to preserve
stability if we restrict the choice of polarisation on $X$ in some
convenient way.

\begin{dfn}
Let $c$ be a positive integer, and consider the set
$$ \Xi(c)=\{ \xi\in{\rm Div}(X)\ |\ -4c\leq\xi^2<0 \ {\rm and}\ \xi{\rm
mod}2=0 \} $$
Let $W^\xi$ be the intersection of the hyperplane $\xi^\perp$ with the
ample cone of $X$. A polarisation $\ell$ is said to be $c$-{\em suitable}
if $\ell\notin W^\xi$ and sign$(\ell\bullet\xi)=$sign$(f\bullet\xi)$
for all $\xi\in\Xi(c)$.
\end{dfn}

It is easy to see that suitable polarizations exist for every $c$.
The following important result is due to Friedman, Morgan \&
Witten \cite{FMW2}:

\begin{thm} \label{rf}
Let $E$ be a torsion-free sheaf on $X$ with $ch(E)=(r,0,-k)$.
If $E$ is $\mu$-semistable with respect to a $k$-suitable
polarisation $\ell$, then $E$ is generically fibrewise semistable.
\end{thm}

\subsection{Preservation of stability}

It follows from Theorem \ref{rf} that if $E$ is a torsion-free sheaf on
$X$ with $ch(E)=(r,0,-k)$, which is $\mu$-semistable with respect
to a $k$-suitable polarisation $\ell$, then $L=\Phi^1(E)$ is
a torsion sheaf of pure dimension one. Let us now consider its
stability (in the sense of Simpson).

\begin{prop} \label{oneway}
Let $E$ be a torsion-free sheaf on $X$ which is $\mu$-semistable with
respect to a $k$-suitable polarisation $\ell$.
Then $L=\Phi^1(E)$ is $p$-semistable with respect to $\hat{\ell}$.
\end{prop}

\begin{proof}
Suppose that
$$ 0\lseta M\lseta L\lseta N\lseta0 $$
is a destabilizing sequence for $L$. We can assume $M$ is
semistable and $N$ has pure dimension 1. Now $p(M,n)=n+\alpha$,
for some $\alpha>0$. Then $\Hom(M,\pp_x)=0$ for all $x\in X$,
since $p(\pp_x,n)=n$. Thus $M$ is $\hat\Phi$\wit0 by Lemma
\ref{tom's}. Arguing as in Proposition \ref{one}, we have
$c_1(\hat\Phi^0(M))=\alpha f$ hence $\mu(\hat\Phi^0(M))>0$. So
$\hat\Phi^0(M)$ will destabilize $E=\hat\Phi^0(L)$ unless
$r(\hat\Phi^0(M))=r(E)$. For this to be the case we must have
$c_1(M)\cdot\hat f=c_1(L)\cdot\hat f$. Then $N$ is supported on a
sum of fibres. But since $N$ must be $\hat\Phi$\wit0 we see that
it has non-negative degree on these fibres and so $\alpha=0$, thus
contradicting the assumption on $M$.
\end{proof}

If we examine the proof more carefully, we can also see that the
assumption that $E$ is $\mu$-stable implies that $L$ is $p$-stable unless
its support decomposes as $D+D'$, where $D$ is the sum of fibres.
When this happens, we have a subsheaf $K$ of $L$ supported on
$D$ with degree 0. Assume that $D$ is the maximal such effective
subdivisor of $S$. Let $Q=L/K$. Then $Q$ is $\hat\Phi$\wit0 and $K$ is
$\hat\Phi$\wit1. The resulting sequence after transforming with
$\hat\Phi$ is then $E\to E^{**}\to\oo_Z$ for some zero dimensional
subscheme $Z$. Conversely, applying $\Phi$ to this sequence shows that
the support of $L$ contains fibres. We have therefore established:

\begin{prop} \label{Estable}
Suppose $E$ is $\mu$-semistable with respect to a $k$-suitable
polarisation $\ell$ with $c_1(E)=0$. If $E$ is $\mu$-stable and
locally-free then $L$ is $p$-stable. Moreover, $E$ is locally-free
if and only if $L$ is destabilized by a sheaf supported on a
fibre.
\end{prop}

We can now consider the opposite question: if we assume $L$ is
$p$-semistable what can we say about its transform?

\begin{prop}\label{Lwit}
Suppose $L$ is a $p$-stable sheaf on $J_X$ with Chern character
$(0,k\hat{f}+r\hat{\sigma},0)$, where $r,k>1$.
Then $L$ is $\hat\Phi$\wit0 and $\hat{\Phi}^0(L)$ is locally-free,
$\mu$-stable with respect to a $k$-suitable polarisation and such
that $ch(\hat{\Phi}^0(L))=(r,0,-k)$.
\end{prop}

\begin{proof}
The first statement follows from lemma \ref{tom's}, which requires
us to show that $\Hom(L,\pp_x)=0$ for all $x$. However, $p(\pp_x,n)=n$
and any map $L\to\pp_x$ would contradict the stability of $L$.

Now suppose that $\hat\Phi^0(L)$ is not $\mu$-stable and let $A$
be the destabilizing subsheaf, so that $\mu(A)\geq\mu(E)=0$
Moreover, we may assume that $A$ is $\mu$-stable, thus $A$ is
$\Phi$\wit1.

We may also assume that the quotient $B=E/F$ is $\mu$-stable.
Then both $A$ and $B$  are $\Phi$\wit1. Since their transforms
must have zero rank, $A$ and $B$ both have zero fibre degree
and so $p(\Phi^1(A),n)\geq n$ and this contradicts the stability
of $L$.

The fact that $\hat{\Phi}^0(L)$ is locally-free follows from the
last part of Proposition \ref{Estable}.
\end{proof}

A similar argument shows that if $L$ is $p$-semistable and
$\hat\Phi$\wit0 then $\hat\Phi^0(L)$ must be $\mu$-semistable.
Note, however, that it the semistability of $L$ does not imply
that $L$ is $\hat\Phi$\wit0. Examining the proof above we see that
this happens precisely when $L$ is destabilized by mapping to a
sheaf supported on a fibre. In particular, such sheaves are
S-equivalent to sheaves which are $\hat\Phi$\wit0 but whose
transforms are not locally-free.

We summarize the results of this section in the following theorem
\begin{thm} \label{summary}
Suppose $E$ is a torsion-free sheaf with rank $r$, $c_1(E)=0$ and
$c_2(E)=k$ on a relatively minimal elliptic surface $X$ over $B$.
Let $\omega$ be a $k$-suitable polarisation on $X$. Then
\begin{enumerate}
\item $E$ is a $\mu$-stable locally free sheaf if and only if its transform
is a $p$-stable torsion sheaf supported on a divisor in
$|k\hat{f}+r\hat{\sigma}|$, where $\hat{\sigma}$ is a section of
$J_X$ and $\hat{f}$ is a fibre class.
\item $E$ is $\mu$-stable properly torsion-free if and only if its
transform is a $p$-semistable torsion sheaf supported on a
reducible divisor in $|k\hat{f}+r\hat{\sigma}|$ which is
destabilized only by a sheaf supported on a fibre.
\item $E$ is properly $\mu$-semistable if and only if its transform is
a $p$-semistable torsion sheaf supported on a reducible divisor in
$|k\hat{f}+r\hat{\sigma}|$.
\end{enumerate}
\end{thm}

One can also give criteria for the Gieseker stability of $E$ in
terms of destabilizing properties of $L$, but these seem to be
less useful.

\begin{rem}
Recently, Hern\'andez Ruip\'erez and Mu\~noz Porras
\cite{HRMP} and Yoshioka \cite{Y} have independently obtained
stability results close to those presented in this Section.
\end{rem}

\begin{rem}
It follows from the first item in Theorem \ref{summary} that the Fourier-Mukai
functor $\Phi$ induces a bijective map from $\cM_X(r,k)$ onto
$\cS_{J_X}(0,k\hat{f}+r\hat{\sigma},0)$, the Simpson moduli space of
p-stable torsion sheaves on $J_X$ with the given Chern classes. As we will
observe in Section \ref{hyper} below, this map is also a hyperk\"ahler
isometry when $X$ is an elliptic K3 or abelian surface (i.e. when $X$ is
hyperk\"ahler).
\end{rem}

We can see that the geometry of the spectral data
is easily linked to the sheaf theoretic properties of the original
sheaf. This should make it very easy in practice to use the spectral
data to analyse properties of the specimen sheaf. We shall see an
example of this in the subsequent Sections where we use the spectral
data to explore the fibration structure on the moduli spaces.


\section{The Fibration Structure} \label{fibre}

From now on we assume that $X$ is either an elliptic $K3$ surface
with a section or a product of elliptic curves. We have seen that
there is a natural map $\Pi$ from the moduli space of $\mu$-stable torsion-free
sheaves $\cM^{\rm TF}_X(r,k)$ to the set $\mathcal{B}$ of spectral curves. In the
case of the K3 surface the base $\projb$ is just the linear system $|k\hat f+r\sigma|$
while for the abelian surface it can be expressed as the total space of the projectivized bundle
$\mathbb{P}\mathcal{F}'(\oo(kf+r\hat\sigma))$, where $\mathcal{F}'$ is the Mukai
transform $D(\hat V\times W)\to D(V\times\hat W)$ (this can be factored as
$\Phi\hat\Upsilon=\hat\Psi\Xi$ from (\ref{commdiag})). The fibre of $\Pi$ over
$S$ is given by suitable subspaces of Jac$_{g(S)-1}(S)$. The
hypothesis on $X$ implies that for $\cM$ to be non-empty we must have
both $r$ and $k$ at least 2.

\begin{thm}
The map $ \cM_X(r,k) \stackrel{\Pi}{\seta} \projb $ is a surjective
map of varieties.
\end{thm}

\begin{proof}
We must first show that $\Pi$ is a well defined map of varieties. To
see this we use the argument given by Friedman and Morgan in
\cite[Chapter VII, Thm 1.14]{FM}. The key idea is to observe that
$\Pi$ coincides (locally in the \'etale or complex topology) with the
projection map from the universal sheaf corresponding to the relative
Picard scheme of degree $g-1$ line bundles on families of genus $g$
curves. The universal sheaf only exists locally in these topologies
but this is enough to show that $\Pi$ is holomorphic.

Given a spectral curve $S$, observe that its structure sheaf $\oo_S$ is
{$\hat{\Phi}$\wit0} since in the short exact sequence
$0\seta\Lambda^{-1}\seta\oo\seta\oo_S\seta0$, $\Lambda^{-1}$ is
$\hat{\Phi}$\wit1 and $\oo$ is $\hat{\Phi}$\wit0.

We aim to construct a $g(S)-1$ degree line bundle $L$ over the curve $S$,
which is stable as a torsion sheaf on $J_X$. That is, we need to choose $g(S)-1$
points on $S$ and guarantee that the restriction of $L$ to any proper component
$S_i$ of $S$ satisfies deg$(L|_{S_i})>0$. We do this by choosing a $g(S_i)$
points on each component $S_i$.  This can always be done because if $g(S_i)>0$
then $S_i$ intersects $S-S_i$ at least twice since $k$ and $r$ are both at least
2 and so $g(S-S_i)<g(S)-g(S_i)$. Without loss of generality, we assume that the
set we have just chosen $Z\subset S$ consists of  distinct points away from the
singularities of $S$. Therefore,
$$ \ext^1(\oo_Z,\oo_S) = \bigoplus_{z\in Z} \ext^1(\oo_z,\oo_S) =
\bigoplus_{z\in Z} T_z S $$
We pick a class in $\ext^1(\oo_Z,\oo_S)$ which is non-zero on each factor. This
defines a torsion sheaf $L$ on $J_X$ given by this extension class. Then $L$ is
stable and locally free on its support by the choice of $Z$.

By Proposition \ref{Lwit}, $L$ is $\hat{\Phi}$\wit0 and
$E=\hat{\Phi}^0(L)$ is $\mu$-stable vector
bundle such that $\Pi(E)=S$.
\end{proof}

From this proof we can also see that any $E\in\cM$ can be written as an
extension
$$ 0 \lseta \hat{\Phi}^0(\oo_S) \lseta E \lseta \hat{\Phi}^0(\oo_Z) \lseta 0. $$
We shall make use of this in the next Section.  Note that such
representation is not unique.

Such a fibration structure will exist for the moduli spaces over any
elliptic surface with a section but this surjectivity result may not
hold. For the rank 2 case see \cite[Thm.~37]{F} or \cite[Thm.~1.14]{FM}.
It is also the case that the fibres will not exist in the middle
dimension which is the situation we wish to consider in the next Section.


\section{The Hyperk\"ahler Structure} \label{hyper}

When $X$ is K3 or a torus, it is well known that the moduli spaces of sheaves are hyperk\"ahler.
To see this is, note that by a result of Mukai's \cite{M}, for each complex structure $I$ on $X$
we obtain a complex symplectic structure $\Omega_I$ on $\cM_X$. Then by Yau's proof of the
Calabi conjecture to obtain the full hyperk\"ahler structure on $\cM_X$.

In fact, the complex symplectic structures arise in a very natural way:
\[ T_{[E]}\cM\times T_{[E]}\cM\cong\ext^1(E,E)\times\ext^1(E,E)
   \stackrel{\cup}{\lseta}\ext^2(E,E)\cong\cpx,\]
where the $E$ is an $\oo_X$-module with respect to $I$.
But we can express this even more simply by observing that
$\ext^i(E,E)=\Hom_{D(X)}(E,E[i])$.
Then the cup product $\cup$ is just composition of maps in the derived category.
Since $\Phi$ is a functor, it must preserve these and so we see that $\Phi$ induces
a complex symplectomorphism:
$$ \cM_X(r,k) \longrightarrow \cS_{J_X}(0,k\hat{f}+r\hat{\sigma},0) $$
Since this happens for each complex structure, the moduli map
induced by $\Phi$ is actually a hyperk\"ahler isometry.

We aim to show now that the fibration structure we have defined in the last
Section on $\cM$ has Lagrangian fibres with respect to this complex symplectic
structure.

\begin{prop}
If $t_1$ and $t_2$ are two tangent vectors to the fibre of $\Pi$ then
$\Omega_I(t_1,t_2)=0$.
\end{prop}

\begin{proof}
By continuity it suffices to prove this when $t_i$ are defined over a point $E$
given by an extension:
\[ 0 \lseta \hat{\Phi}(\oo_S) \lseta E \lseta \hat{\Phi}(\oo_Z) \lseta 0. \]
where $Z$ consists of discrete points and $S$ is smooth. Deformations of
$E$ arising from the fibre of $\Pi$ are determined by deformations of $Z$
along $S$. Then $\ext^2(E,E)\cong\ext^2(L,L)$ is generated by a non-zero
vector in $\ext^2(\oo_Z,\oo_Z)$ and the fibre $\Pi_S$ over $S$ has tangent
space given by $\bigoplus_{z\in Z}\langle\lambda_z\rangle,$ where $\lambda_z$
generates the tangent space to $S$ at $z$. But $\lambda_z\cup\lambda_z=0$ in
$\ext^2(\oo_Z,\oo_Z)$ for each $z$ and so $\Omega_I(t_1,t_2)=0$ as required.
\end{proof}

In the case of a product abelian surface we have natural hyperk\"ahler
actions of the torus and its dual via translations and twisting by
flat line bundles. These act naturally on the fibration structure and
the resulting quotient is also a Lagrangian fibration.

We can say a little more in the particular case of the product of two
elliptic curves. For $\torus=V\times W$, observe that $\cM_\torus(r,k)$
has two such fibration structures:
\begin{equation} \label{dbl.fib} \xymatrix{
& \cM_\torus(r,k) \ar[dl] \ar[dr] & \\
\projb & & \hat{\projb}
} \end{equation}
where $\hat\projb=\mathbb{P}\hat\mathcal{F}(\oo(k\hat f+r\sigma))$.
Note that the fibres are Lagrangian with respect to the same complex
symplectic structure and the bases are biholomorphic since the
underlying vector bundles can be canonically identified by pulling
back along the isomorphism $V\times \hat W\to \hat V\times W$.

In summary we have proved:
\begin{thm}
If X is an elliptic K3 surface with a section or a product of elliptic curves,
then the moduli space of instantons admits a fibration structure over a
compact base (which is a projective space or a projective bundle over an
abelian surface, respectively) and the fibres are Lagrangian with respect
to the natural complex symplectic structure.
\end{thm}

\section*{Acknowledgements}

The first named author would like to thank Yale University and the
University of Pennsylvania for their hospitality.
The second named author would like to thank Yale University and the
Engineering and Physical Sciences Research Council of the UK for their
support. We would also like to thank Tom Bridgeland, Robert Friedman
and the referee for several useful comments.

\end{document}